\documentstyle[12pt]{amsart}

\begin{document}

\title
[Absolutely continuous spectrum of perturbed Stark operators]
{Absolutely continuous spectrum\\of perturbed Stark operators}
\author{Alexander Kiselev}
\address{\hskip-\parindent
Mathematical Sciences Research Institute \\
1000 Centennial Drive \\
Berkeley CA 94720} 

\begin{abstract}
We prove new results on the  stability of the absolutely continuous spectrum for perturbed Stark operators with decaying or satisfying certain smoothness assumption  perturbation.  We show that  the absolutely continuous spectrum of the Stark operator is stable if the perturbing potential decays at the rate $(1+x)^{-\frac{1}{3}-\epsilon}$ or if it is continuously differentiable with derivative from  the H\"older space $C_{\alpha}(R),$ with any $\alpha>0.$
\end{abstract}

\maketitle

\addtocounter{section}{-1}
\section{Introduction}

In this paper, we study the stability of the absolutely 
continuous spectrum of one-dimensional Stark operators
under various classes of perturbations. Stark Schr\"odinger operators describe behavior of the charged particle in the constant electric field. The absolutely continuous spectrum is a manifestation of the fact that the particle described by the operator propagates to infinity at a rather fast rate (see, e.g. \cite{AS}). It is therefore interesting to describe the classes of perturbations which preserve the absolutely continuous spectrum of the Stark operators. In the first part
of this work, we study perturbations of Stark operators by decaying potetnials.
This part is inspired by the recent work of Naboko and Pushnitski \cite{NaPu}.
The general picture that we prove is very similar to the 
case of perturbations of free Schr\"odinger operators \cite{Kis1}. 
In accordance with physical intuition, however, the absolutely continuous spectrum is stable under stronger perturbations than in the free case. 
If in the free case the short range potentials preserving purely absolutely continuous spectrum of the free operator are given by condition (on the power scale)
$|q(x)| \leq C(1+|x|)^{-1-\epsilon},$ in the Stark operator case the corresponding condition reads $|q(x)| \leq C(1+|x|)^{-\frac{1}{2}-\epsilon}.$ If $\epsilon$ is allowed to be zero in the above bounds, imbedded eigenvalues may occur in both cases (see, e.g. \cite{NaPu}, \cite{WN}). Moreover, in both cases if we allow potential to decay slower by an arbitrary function growing to infinity, very rich singular spectrum, such as a dense set of  eigenvalues, may occur (see \cite{Na} for the free case and \cite{NaPu} for the Stark case for precise formulation and proofs of these results). The first part of this work draws the paprallel further, showing that the absolutely continuous spectrum of Stark operators is preserved under perturbations satisfying  $|q(x)| \leq C(1+|x|)^{-\frac{1}{3} -\epsilon},$ in particular even in the regimes where a dense 
set of eignevalues occurs; hence in such cases these eigenvalues are genuinely imbedded. Similar results for the free case were proven in \cite{Kis1}, \cite{Kis2}. Our main strategy of the proof here is similar to that in \cite{Kis1} and \cite{Kis2}: we study the asymptotics of the generalized eigenfunctions and then apply Gilbert-Pearson theory \cite{GiPe} to derive spectral consequences.

While the main new tool we introduce in our treatment of Stark operators is the same as in the free case, namely  the a.e. convergence of the Fourier-type integral operators, there are some major differences. First of all, the spectral parameter enters the final equations that we study in a different way and this makes analysis more complicated. Secondly, we employ a different method to analyze the asymptotics. Instead of Harris-Lutz asymptotic method we study appropriate Pr\"ufer transform variables, simplifying the overall consideration. 

 In the second part of the work we 
discuss perturbations by potentials having some additional smoothness properties, but without decay. It turns out that for Stark operators the effects of decay or of additional smoothness of potential on the spectral properties are somewhat similar. 
It was known for a long time that if a potential perturbing Stark operator has two bounded derivatives the spectrum remains purely absolutely continuous (actually, 
certain growth of derivatives is also allowed, see Section 3 for details or Walter \cite{Wa} for the original result). We note that the results similar to Walter's on the preservation on absolutely continuous spectrum were also obtained in \cite{Ben} by applying different type of technique (Mourre method instead of studying asymptotics of solutions).
  On the other hand, if the perturbing potential is a sequence of derivatives of $\delta$ functions in integer points on $R$
with certain couplings,   the spectrum may turn pure point \cite{AEL}, \cite{Ex}.
 In some sense, the $\delta'$ interaction is the most singular and least
``differentiable" among all available natural perturbations of one-dimensional Schr\"odinger operators \cite{KiSi}. Hence we have very different spectral properties on the very opposite sides of the smoothness scale.  This work closes the part of the gap.  We improve the well-known results of Walter
\cite{Wa}  concerning the minimal smoothness required for the preservation of the absolutely continuous spectrum and show that in fact  existence and minimal smoothness of the first derivative is sufficient to imply absolute continuity of the spectrum.

\section{Decaying perturbations}

Consider a self-adjoint operator $H_{q}$ defined by the differential 
expression 
\[ H_{q}u= -u''-xu+q(x)u \]
on the $L^{2}(-\infty,\infty).$  Let us introduce some notation. For the function $f \in L^{2}$ we denote by $\Phi f$ its Fourier transform:
\[ \Phi f(k)= L^{2}-\lim_{N\rightarrow \infty} \int\limits_{-N}^{N} \exp(ikx)f(x)\,dx. \] 
For locally integrable function $g$ we denote by $M^{+}g$ the function
\[ M^{+}g(x) = \sup_{1>h>0} \frac{1}{2h} \int\limits_{0}^{h} |g(x+t)+g(x-t)|\,dt. \]
We denote by ${\cal M}^{+}f$ the set where $M^{+}f$ is finite. By the general results on the maximal functions (see, e.g. \cite{Rud}), it is easy to conclude that the complement of the set ${\cal M}^{+}$ has Lebesgue measure zero. 
We will prove the following theorem: \\
\bf Theorem 1.1 \it Suppose that the potential $q(x)$ satisfies $|q(x)| \leq C(1+x)^{-\frac{1}{3}-\epsilon},$ $\epsilon >0.$
Then the absolutely continuous spectrum of multiplicity one fills the whole real axis. Imbedded singular spectrum may occur, but only on the complement of the set 
\[ S= {\cal M}^{+}\left(\Phi \left(\exp \left(ix^{3}-i\int_{0}^{x^{3}}q(ct^{\frac{2}{3}})ct^{-\frac{2}{3}}\,dt\right)q(c x^{2})x^{\frac{1}{6}}\right)\right), \]
where $c= (\frac{3}{2})^{\frac{2}{3}}.$

\rm We will prove Theorem 1.1 by studying the asymptotics of the  solutions of the equation 
\begin{equation}
 -u'' -xu+q(x)u=\lambda u .
\end{equation}
We will then apply the subordinacy theory developed by Gilbert and Pearson \cite{GiPe} to derive spectral consequences. For future reference, we summarize here one of the results of the subordinacy theory for Schr\"odinger operators defined on the whole axis  due to Gilbert \cite{Gi} that we will use. 

Let us call the solution $u_{1}$ of the equation (1) subordinate on the right if for any different solution $u_{2}$ the limit 
\[ \lim_{N \rightarrow +\infty} \frac{\|u_{1}\|_{L^{2}(0,N)}}{\|u_{2}\|_{L^{2}(0,N)}} \]
is equal to zero.  Subordinacy on the left is defined similarly. The set of the energies at which there exists a  solution subordinate both on the right and on the left constitutes the essential support of the absolutely continuous part of the spectral measure  of Schr\"odinger operator. Moreover, if for every $\lambda$ from this set there exists a solution subordinate either on the right or on the left, the absolutely continuous spectrum has multiplicity one. 

Since in our case the total potential grows to plus infinity  on the negative  semi-axis, it is clear that for every $\lambda$ there is a solution subordinate on the left and hence we need only to study the behavior of solutions on the positive semi-axis.

In studying Schr\"odinger operators with  sufficiently smooth potentials going to $-\infty$ as $x \rightarrow \infty$ one of the useful tools is Liouville transform (see, e.g.
\cite{Har}). Suppose that  the equation is given by 
\[ -y''+(v(x)+q(x))y= \lambda y, \]
where  $v$ is two times differentiable and goes to $-\infty$ as $x \rightarrow \infty$ sufficiently fast (specifically, $\int_{1}^{\xi} dx/|v(x)|^{\frac{1}{2}}$ diverges as $\xi \rightarrow \infty$) and $q$ is a perturbation. The Lioville transformation is given
by
\[  \xi(x) = \int\limits_{0}^{x} \sqrt{|v(s)|} \,ds, \,\,\,\,\,
\rho(\xi) = |v(x(\xi))|^{\frac{1}{4}}y(x(\xi)). \]
The function $\rho$ satisfies the equation
\begin{equation}
-\rho''+ Q(\xi) \rho = \rho,
\end{equation}
where
\[ Q(\xi)= \frac{5 |v'(x(\xi))|^{2}}{16 |v(x(\xi))|^{3}} - \frac{v''(x(\xi))}{4|v(x(\xi))|^{2}}+\frac{q(x(\xi))-\lambda}{|v(x(\xi))|}. \]
In some cases this new equation does not contain infinite potential and 
is simpler to deal with.

In our case $v(x)=-x$ and we  bring the equation (1) to a convenient form using the following explicit Liouville transformation (see also \cite{NaPu}):
$\xi = \frac{2}{3}x^{\frac{3}{2}},$ $\omega (\xi) = x(\xi)^{\frac{1}{4}}u(x(\xi)).$ Then the function $\omega$ 
solves the  Schr\"odinger equation
\begin{equation}
 -\omega'' + \left( \frac{5}{36\xi^{2}} + \frac{-\lambda + q(c\xi^{\frac{2}{3}})}{c\xi^{\frac{2}{3}}}\right)\omega = \omega 
\end{equation}
(we remind that $c=(\frac{3}{2})^{\frac{2}{3}}$).
We introduce the short-hand notation
\[ b(\xi)= \frac{5}{36\xi^{2}}-\frac{cq(\xi^{\frac{2}{3}})}{c\xi^{\frac{2}{3}}} \]
and $V(\xi, \lambda)$ for the total potential in the equation (3).
Let us further apply Pr\"ufer transformation to the equation for $\omega,$ setting for each $\lambda$
\[ \omega(\xi, \lambda) = R(\xi, \lambda) \sin (\theta(\xi, \lambda)) \]
\[ \omega'(\xi, \lambda) = R(\xi, \lambda) \cos (\theta(\xi, \lambda)). \]
The equations for $R$ and $\theta$ are as follows: 
\begin{equation}
 (\log R(\xi, \lambda))' = \frac{1}{2}V(\xi, \lambda)\sin (2\theta(\xi, \lambda)) 
\end{equation}
\begin{equation}
\theta'(\xi, \lambda) = 1- \frac{1}{2}V(\xi, \lambda)(1-\cos (2\theta(\xi, \lambda))). 
\end{equation}
Our goal now is to study the asymptotics of  solutions of (3).
In particular, we would like to establish convergence  of the integral
\begin{equation}
 \int\limits_{0}^{N} \left( \frac{5}{36\xi^{2}} + \frac{-\lambda + q(c\xi^{\frac{2}{3}})}{c\xi^{\frac{2}{3}}}\right)\sin (2\theta(\xi, \lambda))\, d \xi 
\end{equation}
as $N \rightarrow \infty$ for a.e.$\lambda.$ This goal is motivated by the following  \\
\bf Lemma 1.2. \it Suppose that for some value of $\lambda$ the integral (6) converges. Then for this value of $\lambda,$ there is no subordinate solution of the original equation (1). Moreover, if the integral (6) converges for a.e.~$\lambda,$ the absolutely continuous part of the spectral measure of the operator $H_{q}$ fills the whole real axis.  \\
\bf Proof. \rm If for a given value of $\lambda$ the integral (6) converges,
it follows from  (4) and (5) that all solutions of the equation (3) are bounded  and moreover there are two linearly independent solutions with the following asymptotics as $\xi \rightarrow +\infty:$
\[ \omega_{1}(\xi, \lambda) = \sin (\xi + r_{1}(\xi, \lambda))(1+o(1)), \]
\[ \omega_{2}(\xi, \lambda)= \cos(\xi + r_{2}(\xi, \lambda))(1+o(1)), \]
where $r_{i}'(\xi, \lambda) \leq C(1+ \xi)^{-\frac{2}{3}}.$
Going back to the original equation (1), we infer that for a.e. $\lambda$
this equation has only solutions $u_{\alpha}(x, \lambda)$ (where $\alpha$ is some paprametrization of solutions, say by boundary condition at zero) with the following asymptotical behavior as $x \rightarrow \infty:$
\[ u_{\alpha}( x, \lambda) = x^{-\frac{1}{4}}\sin (\frac{2}{3}x^{\frac{3}{2}}+ f_{\alpha}(x, \lambda))(1+o(1)), \]
where $|f_{\alpha}'(x, \lambda)|\leq C(1+x)^{-\frac{1}{2}}$ uniformly over $\alpha.$ 
For any $u_{\alpha}$ we then find
\[ \int\limits_{0}^{N} |u_{\alpha}(x)|^{2}\,dx = \int\limits_{0}^{N}
x^{-\frac{1}{2}}(\sin(\frac{2}{3}x^{\frac{3}{2}}+ f_{\alpha}(x))^{2}\, dx \left( 1+o(1)\right) = \]
\[ = \left( N^{\frac{1}{2}} + \frac{1}{2}\int\limits_{0}^{N}x^{-\frac{1}{2}}
\cos (\frac{4}{3}x^{\frac{3}{2}}+2 f_{\alpha}(x, \lambda))\, dx \right)\left(1+o(1)\right)= N^{\frac{1}{2}}\left( 1+o(1)\right). \]
The last equality follows from integrating by parts the integral in the previous expression. Hence, we obtain that for a.e.$\lambda$ there are no subordinate solutions: as $N \rightarrow \infty,$ the $L^{2}$ norm grows at the same rate for all solutions. The last claim of the lemma follows  by direct  application of subordinacy theory. $\Box$ \\
\it Remark. \rm 1. Note that in particular we obtained that for a.e. $\lambda$ all solutions of the equation (1) are bounded (and even power-decaying). However the results deriving absolute continuity  of the spectrum from the boundedness of solutions (see, e.g., \cite{Sto}, \cite{Si2}) are not applicable here because the potential goes to $-\infty$ on $R^{+}.$ Instead, one has to apply directly Gilbert-Pearson theory. \\
2. Note that while for a.e.$\lambda$ there is no subordinate solution of equation (1) and hence the absolutely continuous spectrum fills the whole axis, there may be a complementary set of measure zero where the singular part of the spectral measure might be supported. As examples show \cite{NaPu}, the imbedded spectrum may even be dense. The a.e. convergence of Fourier type integral operators  allows to prove a.e. convergence of the integral (6) while permitting for exceptional measure zero set where convergence may fail. \\

We begin analyzing the integral (6) with the following simple \\
\bf Lemma 1.3. \it Suppose that the function $h(\xi)$ satisfies 
\[ h(\xi)= \xi + g(\xi), \]
where $|g'(\xi)|\leq C\xi^{-\frac{2}{3}}.$ 
Then the integrals
\[ \int\limits_{1}^{N} \xi^{-\frac{2}{3}} \exp (\pm i h(\xi))\, d\xi \]
converge and , moreover, 
\[ \int\limits_{N}^{\infty}\xi^{-\frac{2}{3}} \exp (\pm i h(\xi))\, d\xi= O(N^{-\frac{1}{3}}). \]
\bf Proof. \rm 
\[ \int\limits_{1}^{N} \xi^{-\frac{2}{3}} \exp (\pm i h(\xi))\, d\xi = N^{-\frac{2}{3}}\int\limits_{1}^{N} \exp (\pm i h(\xi))\,d\xi +\]
\[ + \frac{2}{3}\int\limits_{1}^{N}\xi^{-\frac{5}{3}}\int\limits_{1}^{\xi}
\exp (\pm i h(\eta))\, d\eta d\xi. \]
 Pick $\xi_{0}$ so that 
for $\xi>\xi_{0}$ $|g'(\xi)|<\frac{1}{2}$ ( we can do it by the assumption of Lemma). Then 
\[ \int\limits_{\xi_{0}}^{N} \exp(\pm i h(\xi))\,d\xi=
\int\limits_{h(\xi_{0})}^{h(\xi)}\exp(\pm i h)\frac{1}{1+g'(\xi(h))}\, dh, \]
where $|g'(\xi(h)))| \leq C_{1}h^{-\frac{2}{3}}.$
Expanding the fraction in the last integral, it we find that 
\[ \int\limits_{\xi_{0}}^{N} \exp(\pm  i h(\xi))\,d\xi=O(N^{\frac{1}{3}}). \]
Hence we see that the integral in question converges.
It is also straightforward to establish the last estimate of the lemma. 
$\Box$

From Lemma 1.3 and (5) it follows that the integrals 
\[ \int\limits_{0}^{N} \xi^{-\frac{2}{3}} \exp(\pm 2i \theta(\xi, \lambda)) \,d\xi \]
converge for every value of $\lambda.$
Therefore, we see  that to establish convergence  of the integral (6)
it suffices to study  the integral
\begin{equation}
 \int\limits_{1}^{N} \left( \frac{5}{36\xi^{2}}-q(c\xi^{\frac{2}{3}})(c\xi)^{-\frac{2}{3}}\right)\sin(2\theta(\xi, \lambda)\, d\xi). 
\end{equation}
(Of course, for the fast decaying term  the convergence issue is trivial, but  we incorporate it into the expression for future computational convenience.) 

 We will need two  lemmata on the a.e. convergence of Fourier-type integrals. The first lemma we need is exactly Lemma 1.3 from \cite{Kis1}.
We refer to that paper for a proof. \\
\bf Lemma 1.4. \it Consider the function $f \in L^{2}(R).$ Then for every $\lambda_{0} \in {\cal M}^{+}(\Phi(f))$ we have
\[ \int\limits_{-N}^{N} f(x) \exp(i\lambda_{0}x)\,dx = O(\log N). \]

\rm The second lemma is a variation of 
Lemma 1.4 from \cite{Kis1}. We provide here the statement and a sketch of the proof for the sake of completness. \\
\bf Lemma 1.5. \it Suppose that a function $f(x)$ satisfies $|f(x)| \leq 
C(1+|x|)^{-\alpha},$ where $\alpha > \frac{1}{2}.$ Then for a.e. $\lambda$ the integral 
\[ \int\limits_{0}^{N} \exp(i\lambda x)f(x)\, dx \]
converges as $x \rightarrow \infty$ and moreover for every $\delta >0$ 
\begin{equation}
 \int\limits_{N}^{\infty} \exp(i\lambda x)f(x)\, dx = o(N^{-\alpha+\frac{1}{2}+\delta}) 
\end{equation}
for a.e. $\lambda.$ Explicitely, we can say that (8) holds for every $\lambda \in {\cal M}^{+}(\Phi (f(x)x^{\alpha-\frac{1}{2}-\epsilon}))$ where $\epsilon<\delta.$ \\
\bf Proof. \rm   Consider an apriori estimate 
\[ \int\limits_{N}^{\infty} \exp (i\lambda x) x^{\alpha-\frac{1}{2}-\epsilon}f(x)x^{-\alpha+\frac{1}{2}+\epsilon}\,dx= 
N^{-\alpha+\frac{1}{2}+\epsilon}\int\limits_{0}^{N}\exp(i\lambda x)
(f(x)x^{\alpha-\frac{1}{2}-\epsilon})\,dx+ \]
\[+(\frac{1}{2}-\alpha+\epsilon)\int\limits_{N}^{\infty}x^{-\alpha-\frac{1}{2}+\epsilon}\left(\int\limits_{0}^{x}\exp(i \lambda y)(f(y)y^{\alpha-\frac{1}{2}-\epsilon}\,dy \right)\,dx. \]
By Lemma 1.4, for every $\lambda \in {\cal M}^{+}(\Phi (f(x)x^{\alpha-\frac{1}{2}-\epsilon}))$ we have an estimate
\[ \int\limits_{0}^{N}\exp(i \lambda y)f(y)y^{\alpha-\frac{1}{2}-\epsilon}\,dy =O(\log N) \]
 Therefore for such $\lambda$
we obtain 
\[ \int\limits_{N}^{\infty}\exp(i \lambda x)f(x)\,dx = O(N^{-\alpha+\frac{1}{2}+\epsilon}\log N) = o(N^{-\alpha +\frac{1}{2} +\delta}) \]
since $\delta>\epsilon.$ The lemma is proven. $\Box$ \\
\it Remark. \rm We note that by Zygmund's theorem \cite{Zyg} the Fourier integral converges a.e. for every $f \in L^{p},$ $1 \leq p<2.$ One can avoid using Lemma 1.4 by applying directly this theorem in the proof of Lemma 1.5. However this approach is less direct and also does not provide any explicit information about the convergence set. \\

\noindent \bf Corollary 1.6. \it Suppose that the function $f(\xi)$ satisfies $|f(\xi)| \leq C(1+|\xi|)^{-\alpha-\epsilon},$ where $\alpha \geq \frac{5}{6}$ and $\epsilon$ is an arbitrarily small positive number. Then for a.e.$\lambda$ the integral 
\[ \int\limits_{0}^{N} f(\xi) \exp(i\lambda \xi^{\frac{1}{3}})\, d\xi \]
converges as $N \rightarrow \infty.$ Moreover, for every $\lambda \in {\cal M}^{+}(\Phi (f(\xi^{3})\xi^{3\alpha-\frac{1}{2}}))$ we have
\[ \left| \int\limits_{N}^{\infty} f(\xi)\exp(i\lambda \xi^{\frac{1}{3}})\,d\xi \right| \leq CN^{-\alpha +\frac{5}{6}}. \]
\bf Proof. \rm 
The change of variables $y=\xi^{\frac{1}{3}}$ transforms the integral to the form
\[ \int\limits_{0}^{N^{\frac{1}{3}}} f(y^{3})y^{2}\exp (6i\lambda y)\, dy. \]
Note that 
\[ |f(y^{3})y^{2}| \leq C(1+|y|)^{-3\alpha-3\epsilon+2} \]
by assumption. Applying Lemma 1.5 (and remembering to take into account  the lower limit) we obtain the statements of the Corollary. $\Box$

We prove the final lemma we need for the proof of the theorem. 
Let us introduce some short-hand notation:
\[ \sigma (x, \lambda)= 2\theta(x, \lambda)- 2x + \int\limits_{0}^{x}
\left( q(\xi^{\frac{2}{3}})\xi^{-\frac{2}{3}}-\lambda\xi^{-\frac{2}{3}}\right)\,d\xi+\lambda \int\limits_{0}^{x}\xi^{-\frac{2}{3}}\cos 2\theta(\xi, \lambda)\,d\xi; \]
\[ \gamma(x, \lambda) = 2\theta(x, \lambda) - \sigma(x, \lambda). \] 
From the equation (5) we find
\begin{equation}
\gamma'(x, \lambda)= b(x, \lambda) \cos (\gamma(x, \lambda)+\sigma(x, \lambda)). 
\end{equation}
\bf Lemma 1.7. \it Suppose that the potential $q(x)$ satisfies $|q(x)| \leq C(1+|x|)^{-\frac{1}{3}-\epsilon}.$ Then the integrals 
\[ \int\limits_{0}^{N} b(\xi) \exp (\pm \sigma (\xi, \lambda))\,d\xi \]
converge as $N \rightarrow \infty$ for a.e. $\lambda$ and moreover
for a.e. $\lambda$ 
\begin{equation}
 \int\limits_{N}^{\infty} b(\xi) \exp (\pm \sigma (\xi, \lambda))\,d\xi = O(N^{-\frac{1}{18}}). 
\end{equation}
\bf Proof. \rm We will consider the case of plus sign; the other case is analogous. By Lemma 1.3 and the definitions of $b(\xi)$ and $\sigma(\xi, \lambda),$ it is clear that it suffices to investigate the convergence of the following integral
\begin{equation} \int\limits_{0}^{N} q(c\xi^{\frac{2}{3}})\xi^{-\frac{2}{3}} \exp\left(
2i\xi -i\int_{0}^{\xi}q(\eta^{\frac{2}{3}})\eta^{-\frac{2}{3}}\,d\eta+ 6\lambda i\xi^{\frac{1}{3}} +ic(\lambda)+i\xi^{-\frac{1}{3}}h(\lambda, \xi)\right), 
\end{equation}
where the function $h(\xi, \lambda)$ is bounded for every  $\lambda.$ 
We used Lemma 1.3 to   replace 
\[ \lambda \int\limits_{0}^{\xi} \eta^{-\frac{2}{3}}\cos (2\theta (\eta, \lambda))\, d\eta \]
by \[ c(\lambda) +  \xi^{-\frac{1}{3}}h(\lambda, \xi) \]
in the exponent in (11).
Let us denote 
\begin{equation}
 \tilde{q}(\xi)=q(c\xi^{\frac{2}{3}})\xi^{-\frac{2}{3}}\exp(2i\xi-i\int_{0}^{\xi}q(c\eta^{\frac{2}{3}})c\eta^{-\frac{2}{3}}\,d\eta).
\end{equation}
We can rewrite the integral (11) as 
\[ \int\limits_{0}^{N} \tilde{q}(\xi) \exp (6i\lambda \xi^{\frac{1}{3}} -
i\xi^{-\frac{1}{3}}h(\lambda, \xi))\, d\xi \]
(henceforth we omit the constant factor $\exp (ic(\lambda))$).
Note that as $\xi \rightarrow \infty,$ we have  
\[ \exp(i \xi^{-\frac{1}{3}} h(\lambda, \xi))= 1+ \xi^{-\frac{1}{3}}\tilde{h}(\lambda,\xi), \]
where $\tilde{h}$ is bounded  in $\xi$ for every $\lambda.$
Hence we can write the integral (11) as a sum of two integrals
\[ \int\limits_{0}^{N}\tilde{q}(\xi) \exp (6i\lambda \xi^{\frac{1}{3}})\,d\xi + \int\limits_{0}^{N}\tilde{q}(\xi)\xi^{-\frac{1}{3}}\tilde{h}(\xi, \lambda) \exp (6i\lambda, \xi)\,d\xi. \]
Note that by assumption on the decay of $q$ and (12) we have 
\begin{equation}
 \tilde{q}(\xi) \leq C\xi^{-\frac{8}{9}-\frac{2}{3}\epsilon}. 
\end{equation}
Therefore, the  second integral converges for every $\lambda$ because the integrand is absolutely integrable. It is also easy to estimate the rate of convergence 
for the second integral:
\[ \left|\int\limits_{N}^{\infty}\tilde{q}(\xi)\xi^{-\frac{1}{3}}\tilde{h}(\xi, \lambda) \exp (6i\lambda \xi)\,d \xi \right| \leq C(\lambda)N^{-\frac{2}{9}} \]
for every $\lambda$ because of the assumption on decay of $q.$
To study the behavior of the first integral, we note that because of the estimate (13) we can apply Corollary 1.6.
We obtain that for a.e. $\lambda$ the first integral converges, and, moreover, for a.e. $\lambda$ (explicitely, for every $\lambda \in {\cal M}^{+}(\Phi (\tilde{q}(\xi^{3})\xi^{\frac{13}{6}})$) the estimate on the rate of convergence is given by
\[ \int\limits_{N}^{\infty}\tilde{q}(\xi) \exp (6i\lambda \xi^{\frac{1}{3}})\,d\xi \leq C(\lambda)N^{-\frac{1}{18}}. \]
Combining the two estimates, we see that the lemma is proven.
 $\Box$

We are now in a position to complete the proof of the theorem. The main idea is to take an advantage of certain symmetry of the equations (4), (5) and estimate (6). \\
\bf Proof of Theorem 1.1. \rm  As we already remarked, to complete the proof we only have to study the a.e. convergence of the integral (7). Using (9) we rewrite this integral as
\[ \int\limits_{1}^{\infty} b(\xi) \sin(\gamma (\xi, \lambda)+\sigma(\xi, \lambda))\,d\xi. \]
We will integrate by parts two times, differentiating the terms which contain $\gamma(\xi, \lambda).$ We omit the arguments of $\gamma$ and $\sigma$ in the following computations.
\[ \int\limits_{1}^{N} b(\xi)\sin (\gamma +\sigma)\,d \xi  = \int\limits_{1}^{N} b(\xi)(\sin (\gamma )\cos(\sigma )+ \sin (\sigma)\cos (\gamma))\, d\xi = \]
\[ = -\left( \left. \cos (\gamma) \int\limits_{\xi}^{\infty} b(\eta)
\sin(\sigma)\, d\eta + \sin (\gamma) \int\limits_{\xi}^{\infty} b(\eta) \cos (\sigma)\,d \eta \right) \right|^{N}_{1}+ \]
\[ + \int\limits_{1}^{N} b(\xi)\cos(\gamma+\sigma)\left(\cos (\gamma )\int\limits_{\xi}^{\infty} b(\eta)\cos(\sigma)\,d \eta-\sin (\gamma)\int\limits_{\xi}^{\infty} b(\eta)\sin(\sigma)\,d \eta \right)\,d\xi. \]
We remark that by Lemma 1.7 the conditional integrals are well-defined for a.e. $\lambda.$ Given that, we see that for a.e. $\lambda$ the off-diagonal terms stay finite. We rearrange the integral term and write it as follows:
\[ \int\limits_{1}^{N} (\sin(\gamma))^{2} \left( b(\xi) \sin (\sigma)\int\limits_{\xi}^{\infty} b(\eta)\sin(\sigma)\,d \eta \right)  d \xi - \]
\[- \int\limits_{1}^{N} (\sin(\gamma)\cos (\gamma)) \left( b(\xi) \cos (\sigma)\int\limits_{\xi}^{\infty} b(\eta)\sin(\sigma)\,d \eta +  b(\xi) \sin (\sigma)\int\limits_{\xi}^{\infty} b(\eta)\cos(\sigma)\,d \eta \right) d \xi  + \]
\[ + \int\limits_{1}^{N} \cos(\gamma)^{2} \left( b(\xi) \cos (\sigma)\int\limits_{\xi}^{\infty} b(\eta)\cos(\sigma)\,d \eta  \right) d \xi . \]
Integrating by parts (terms in brackets being integrated), we again  obtain off-diagonal terms which are a.e.$\lambda$ finite and the following integral term:
\[ \int\limits_{1}^{N} \sin (2\gamma) b(\xi) \cos (\gamma+\sigma) \left( \left( \int\limits_{\xi}^{\infty}b(\eta)\sin (\sigma)\,d \eta \right)^{2} - \left( \int\limits_{\xi}^{\infty}b(\eta)\cos (\sigma)\,d \eta \right)^{2}\right)\, d \xi - \]
\[ - \int\limits_{1}^{N} \cos (2\gamma) b(\xi) \cos (\gamma+\sigma)  \left( \int\limits_{\xi}^{\infty}b(\eta)\sin (\sigma)\,d \eta \right)\left(\int\limits_{\xi}^{\infty}b(\eta)\cos (\sigma)\,d \eta \right) \, d \xi. \]
Hence, the integral (7) and therefore (6) is convergent if 
\[ b(\xi) \left( \int\limits_{\xi}^{\infty} b(\eta) \exp (\pm \sigma(\eta, \lambda))\, d \eta \right)^{2} \in L^{1}(0, \infty). \]
Because of the assumption on the decay of $q$ definition of $b(\xi)$ and Lemma 1.7, we find that for a.e. $\lambda \in R$ the expression above is bounded above by 
\[ C(\lambda) \xi^{-\frac{2}{3}}\xi^{-\frac{2}{9}-\frac{\epsilon}{2}}\xi^{-\frac{1}{9}}\leq C(\lambda) \xi^{-1-\frac{\epsilon}{2}} \]
and hence is absolutely integrable for a.e. $\lambda.$  This proves the convergence of the integral (6) for a.e. $\lambda.$ Since from the proof of Lemma 1.7  follows that the estimate (10) holds for every $\lambda \in {\cal M}^{+}(\Phi (\tilde{q}(\xi^{3}))\xi^{-\frac{13}{6}})),$ we also obtain an explicit set to which the singular part of the spectral measure gives zero weight, as stated in the  theorem. $\Box$

\section{The case of sufficiently smooth potential}

In this section, we take up the study of an alternative kind of conditions implying the absolute continuity of the spectrum of Stark operators. 
We begin with the remark that if the perturbation $q(x)$ is twice differetiable, $|q'(x)| \leq C(1+|x|)$ and $|q''(x)| \leq C(1+|x|)^{\alpha},$ $\alpha<\frac{1}{2},$ we can apply the Liouville transformation directly to the 
whole potential $-x+q(x)$, obtaining the equation (2) with absolutely integrable 
potential $Q.$ Then by Lemma 1.2, for every $\lambda$ there is no subordinate solution of the original equation and hence the spectrum is purely absolutely continuous. The situation is more subtle if $q$ does not have two derivatives and hence we cannot apply Liouville transformation directly to the whole potential $-x+q(x).$
Our goal is to show the following improvement of results proven in  \cite{Wa}, \cite{Ben}: \\
\bf Theorem 2.1. \it Suppose that the perturbation $q(x)$ is bounded and differentiable, and the derivative $q'(x)$ is bounded and H\"older continuous with any exponent $\alpha >0,$ i.e.
\[ \sup_{x \in R^{+}} \left(\sup_{y} \frac{|q'(x)-q'(y)|}{|x-y|^{\alpha}}\right)<C. \]
Then the absolutely continuous part of the spectral measure $\rho_{\rm{ac}}^{H_{q}}$ associated with the Stark operator $H_{q}$ fills the whole real axis. \\
\bf Proof. \rm 
The strategy of the proof will be the same as for decaying case. We begin by analyzing equations (4), (5):
\[ (\log R)'(\xi, \lambda)= \frac{1}{2} V(\xi)\sin (2\theta (\xi, \lambda)) \]
\[ \theta'(\xi, \lambda) = 1 -\frac{1}{2}V(\xi)+\frac{1}{2}V(\xi)\cos(2\theta(\xi, \lambda)). \]
We remind that 
\[ V(\xi) =  \frac{5}{36\xi^{2}}+\frac{\lambda - q(c\xi^{\frac{2}{3}})}{c\xi^{\frac{2}{3}}}  \]
To treat the phase equation, we introduce quantities similar to $\gamma$ and $\sigma$ in the previous proof:
\[ \tilde{\sigma}(\xi, \lambda)= 2 \xi - \int\limits_{0}^{\xi} V(\eta)\,d \eta \]
\[ \tilde{\gamma}(\xi, \lambda) = 2\theta (\xi, \lambda) - \tilde{\sigma}(\xi, \lambda). \]
 The computation completely repeating that in the proof of Theorem 1.1
allows us to conclude that $R(\xi, \lambda)$ is bounded for a.e. $\lambda$
if the expression 
\begin{equation}
 V(\xi,\lambda )\left[ \int\limits_{\xi}^{\infty} V(\eta, \lambda) \exp \left( \pm \left( 2i\eta -\int\limits_{0}^{\eta} V(t, \lambda)\, dt \right)\right)\, d\eta \right]^{2} 
\end{equation}
is well-defined and absolutely integrable for a.e. $\lambda.$
Let us now consider this control expression more carefully. 
Consider the  integral 
\[ w_{\pm}(\xi, \lambda)=\int\limits_{\xi}^{\infty} V(\eta, \lambda) \exp \left( \pm \left( 2i\eta -\int\limits_{0}^{\eta} V(t, \lambda)\, dt \right)\right)\, d\eta = \]
\[ = -\int\limits_{\xi}^{\infty} \frac{q(c\eta^{\frac{2}{3}})}{c\eta^{-\frac{2}{3}}} \exp \left( \pm \left( 2i\eta -\int\limits_{0}^{\eta} V(t, \lambda)\, dt \right)\right)\, d\eta \,\,+ \] \[ + \int\limits_{\xi}^{\infty} \left(\frac{5}{36\eta^{2}} +\frac{\lambda}{c\eta^{\frac{2}{3}}} \right) \exp \left( \pm \left( 2i\eta -\int\limits_{0}^{\eta} V(t, \lambda)\, dt \right)\right)\, d\eta. \]
We see that the  second integral on the right-hand side converges for all $\lambda$ and is $O(\xi^{-\frac{1}{3}})$ by  Lemma 1.3 
Hence it remains to study the first integral. We rewrite this integral as follows:
\begin{equation}
 \int\limits_{\xi}^{\infty} q(c\eta^{\frac{2}{3}}) \exp \left( \pm i \left(  -\frac{5}{36\eta}+ \int\limits_{0}^{\eta}q(ct^{\frac{2}{3}})ct^{-\frac{2}{3}}\,dt \right)\right) c\eta^{-\frac{2}{3}}\exp \left( \pm \left(
2i\eta + 6i\lambda \eta^{\frac{1}{3}} \right)\right)\, d\eta. 
\end{equation}
Let us introduce a short-hand notation $g(\eta)$ for the function 
\[ \exp \left( \pm i \left(  -\frac{5}{36\eta}+ \int\limits_{0}^{\eta}q(ct^{\frac{2}{3}})ct^{-\frac{2}{3}}\,dt \right)\right).\]
Note that $g$ is continously differentiable and $|g(\eta)|+|g'(\eta)|$ is bounded.
In (15) we differentiate by parts, differentiating the first two terms and integrating next two. Note that 
\[ \int\limits_{\xi}^{\infty} \eta^{-\frac{2}{3}} \exp \left(\pm \left(
2i\eta +6i\lambda \eta^{\frac{1}{3}} \right)\right) \, d\eta = \frac{1}{2}\xi^{-\frac{2}{3}}\exp \left( \pm 2i\xi \pm 6i\lambda \xi^{\frac{1}{3}} \right) +
O(\xi^{-\frac{5}{3}}). \]
Given this observation, we obtain that the integral (15) is equal to
\[ O(\xi^{-\frac{2}{3}}) + \left( \frac{c}{6}\int\limits_{\xi}^{\infty} 
\left(2\eta^{-\frac{1}{3}}q'(c\eta^{\frac{2}{3}}) \pm i q(c\eta^{\frac{2}{3}})
\eta^{-\frac{2}{3}} \mp i \frac{5}{36c\eta^{2}} \right) \right. \times \]
\[ \times \left. g(\eta)\exp(\pm 2i \eta \pm 6i \lambda \eta^{\frac{1}{3}})
\eta^{-\frac{2}{3}}\left( 1+ O(\eta^{-1})\right) \right) = \]
\[ = O(\xi^{-\frac{1}{3}}) + \frac{c}{6} \int\limits_{\xi}^{\infty}q'(c\eta^{\frac{2}{3}})\eta^{-1}g(\eta)\exp \left( \pm 2i\eta \pm 6i\lambda \eta^{\frac{1}{3}} \right)\, d \eta. \]
By assumption, $q'$ is is H\"older continuous uniformly over $R$ with 
exponent $\alpha.$ Hence for every $h>0,$
\begin{equation}
 q'(c(\eta + h)^{\frac{2}{3}}) - q'(c\eta^{\frac{2}{3}}) = C( \eta, h)\eta^{-\frac{\alpha}{3}}h^{\alpha}, 
\end{equation}
where $C(\eta, h)$ is bounded uniformly for all $\eta,$ $0<h<1.$ In the rest of the proof we will keep notation $C(\eta, h)$ for different functions satisfying the same property. 
Let us integrate by parts the last integral ``$\alpha$" times. That is, let $h>0$ be fixed, such that $\exp(ih) -1 \neq 0.$ Then 
\[ \int\limits_{\xi}^{\infty} q'(\eta^{\frac{2}{3}})g(\eta)\eta^{-1} \exp (\pm 6 i\lambda \eta^{\frac{1}{3}})\exp(\pm 2i\eta) \left( \exp(\pm 2ih)-1 \right)\, d \eta= \]
\[ O(\xi^{-1}) + \int\limits_{\xi}^{\infty} \left( q'((\eta+h)^{\frac{2}{3}})(\eta+h)^{-1}g(\eta+h)\exp( \pm 6i\lambda (\eta+h)^{\frac{1}{3}})- \right. \] \[- \left. q'(\eta^{\frac{2}{3}})\eta^{-1}g(\eta) \exp(\pm 6i\lambda\eta^{\frac{1}{3}}) \right) \exp( \pm 2i\eta)\, d \eta. \]
Note that 
\[ \exp(\pm 6i\lambda(\eta+h)^{\frac{2}{3}})-\exp(\pm 6i\lambda\eta^{\frac{2}{3}}) = \eta^{-\frac{2}{3}}C_{\lambda}(\eta, h). \]
Hence we can rewrite the last integral as
\[ O(\xi^{-\frac{2}{3}}) + \int\limits_{\xi}^{\infty} \left((q'((\eta+h)^{\frac{2}{3}})(\eta+h)^{-1}g(\eta+h)- q'(\eta^{\frac{2}{3}})\eta^{-1}g(\eta)\right)\exp(\pm 2i\eta \pm 6i\lambda \eta^{\frac{1}{3}})\, d \eta.  \]
Since we have
\begin{equation}
|(\eta+h)^{-1}g(\eta+h)-\eta^{-1}g(\eta)| \leq C\eta^{-\frac{2}{3}},
\end{equation}
a simple computation taking into account (16) and (17) shows that the last expression is equal to 
\[  O(\xi ^{-\frac{2}{3}})+\int\limits_{\xi}^{\infty} C(\eta, h)\eta^{-1-\frac{\alpha}{3}}\exp(\pm 2i\eta \pm 6i\lambda \eta^{\frac{1}{3}})\, d \eta.  \]
 By Corollary 1.6, we obtain that the last expression converges for a.e. $\lambda$ and is bounded by $C\xi^{-\frac{1}{3}(\frac{1}{2}+\alpha)+\epsilon}$ for every $\epsilon >0.$
Therefore we have completed an estimation of the integral $w_{\pm}(\xi, \lambda)$ and obtained that for a.e. $\lambda$
\[ w_{\pm}(\xi, \lambda) \leq C(\lambda) (\xi^{-\frac{1}{3}+\epsilon}+
\xi^{\frac{1}{3}(\frac{1}{2}+\alpha)+\epsilon}). \]
Taking into account the fact that $|V(\xi)| \leq C\xi^{-\frac{2}{3}},$ this implies that for a.e. $\lambda$ the control expression (14)
is bounded by $C(\lambda)\xi^{-1-\frac{2\alpha}{3}+\epsilon}$
and hence is absolutely integrable. This allows us to conclude the boundedness of solutions of the equation (3) for a.e. $\lambda.$ Applying Lemma 1.2, we find that for a.e. $\lambda$ there is no subordinate solution of the original equation (1). This completes the proof of Theorem 2.1. $\Box$ \\

\section*{Acknowledgment} I am very grateful to Professor  for stimulating discussions and to Professor Stolz for turning author's attention to the problem.

\end{document}